\documentclass[twocolumn]{IEEEtran}
\usepackage{amsmath, amsthm, graphicx}
\usepackage{subfigure}
\usepackage{amssymb,amsfonts}

 \newcommand{\argmin}{\mathop{\rm argmin}}
 \newcommand{\sign}{\mathop{\rm sign}}

\newtheorem{theorem}{Theorem}
\newtheorem{proposition}{Proposition}

\begin{document}

\title{Adaptive optimal $\ell_\infty$-induced robust stabilization of minimum phase SISO plant under bounded disturbance and coprime factor perturbations}

\author{Victor F. Sokolov\\

Institute of Physics and Mathematics, Komi Scientific Center of Ural Branch of RAS, Syktyvkar , Russia\\
email: vfs-t@yandex.ru}

\maketitle




\begin{abstract}
This paper addresses the problem of optimal robust stabilization of a discrete-time minimum-phase plant in the framework of robust control theory in the $\ell_1$ setup and  under poor a priori information. Coefficients of the transfer function of the plant nominal model with stable zeros are unknown and belong to a known bounded polyhedron in the space of  coefficients. The gains of coprime factor perturbations of the plant and the upper bound  of external disturbance are also unknown. The problem under consideration is to design adaptive controller that minimizes, with the prescribed accuracy, the worst-case asymptotic upper bound of the output. Solution of the problem is based on set-membership estimation of unknown parameters and treating the control criterion as the identification criterion. A hard nonconvex problem of on-line computation of optimal estimates is reduced, under additional nonrestrictive assumption, to a linear-fractional programming via a nonlinear transformation of estimated parameters. Despite the non-identifiability of the unknown parameters,  the proposed adaptive controller guarantees, with the prescribed accuracy, the same optimal asymptotic upper bound of the output of adaptive system as the optimal controller for the plant with known parameters. In addition to the optimality  of adaptive control, the proposed solution provides on-line verification/validation of current estimates and a priori assumptions.
\end{abstract}



\maketitle


\section{Introduction}
\label{sec1}

In this paper, by adaptive we mean control of systems with linear time invariant nominal model, parameters of which are unknown to controller designer and are estimated in closed loop. The estimation of the unknown parameters is typically based on various gradient type algorithms or modifications of the recursive least squares (RLS) algorithm. In the early 1980s, it was realized that adaptive systems with this kind of estimation algorithms can demonstrate unacceptable dynamics or instability for systems with bounded external disturbances and/or uncertainties (unmodelled dynamics)\cite{rohr82}. This motivated the development of the theory of \textit{robust adaptive control}, the main goal of which was to ensure stability of adaptive systems and to derive, whenever possible, some performance guarantees based typically on Lyapunov's methods.  \cite{iosun96}. More difficult problems of \textit{adaptive optimal control} were solved in stochastic settings for systems with random external disturbances and no uncertainties with the use of the gradient and RLS type estimation algorithms \cite{good81,guochen91}. However, there are no extensions of stochastic adaptive optimal control on systems with unmodelled dynamics because even theory of nonadaptive stochastic optimal robust control was not elaborated so far. For systems with unmodelled dynamics and random external disturbances, the stability of adaptive system was  also obtained, but the proof of stability in the mean-square sense was based on undesired bursts of large magnitude due to unmodelled dynamics \cite{rad95}. Note that the theory of robust adaptive control in 1980-s and 1990-s was mainly based on the Lyapunov theory  and was little correlated with the nonadaptive theory of robust control, which was developed in parallel in the same years with the use of other methods like $\mu$-synthesis and $H_\infty$ control  \cite{zhou96}. Therefore, the results on robust stability of adaptive systems with sufficiently small unmodelled dynamics  were rather of qualitative than quantitative nature. More advanced result, based on the small gain theorem, was obtained with the use of gradient type estimation for a special plant with the known bound of the external disturbance and the known gain of unstructured uncertainty \cite{wey94}. However, the gain of the uncertainty and the asymptotic upper bound on the plant output were the same for all parameters from  a priori set of unknown parameters of the plant and, therefore, were very conservative.

The model of bounded external disturbance has generated the set-membership approach to system identification. This approach is based on the assumption of known upper bound on the $\ell_\infty$-norm of the disturbance. For models linear in their parameters, the sets of unfalsified by data parameters  are polytopes in the parameters space and the number of linear inequalities in their description can grows linearly with time. The main problem under consideration was to find limited complexity approximations of these sets via orthotopes, parallelotopes, zonotopes, ellipsoids or others.  Researches on this problem, mostly without any applications to control, were presented in hundreds of papers and several special issues (references \cite{mil96,nor95,nor94,wal90} are only a part of them)  and continue till present time. In recent years, set-membership estimation began to be used in control problems. As an example, parallelotope estimates are used for adaptive model predictive building climate control \cite{tan17}. A data-driven algorithm to approximately compute a minimal robust invariant set (RCI) in the form of the polytope with predefined faces orientation by simultaneously selecting an admissible model and minimizing the size of the RCI is presented for  autonomous vehicle lane-keeping control \cite{chen21}. Polytopic estimates were used in adaptive model predictive control to non-conservatively guarantee recursive feasibility and constraint satisfaction for systems with parametric uncertainty under the assumption of the existence of a robustly stabilizing feedback law and a common Lyapunov function for the stabilized system, for all model parameters in a given prior bounding set \cite{lor19}.

In this paper, the problem of adaptive robust stabilization is considered in the optimal setup. The controlled plant is linear time-invariant (LTI) single-input single-output (SISO) minimum-phase system under bounded disturbance and bounded memory coprime factor perturbations. The coefficients of the transfer function of the nominal system are unknown and belong to a known polytope. The upper bound of the external disturbance and the gains of bounded memory perturbations are also unknown. The control objective is to minimize the worst-case upper limit of the absolute value of the output. The assumption of unknown upper bound of the external disturbance implies, for any time interval, that any coefficients of the transfer function from a priori polytope are not unfalsified by data under sufficiently large norm of the external disturbance. So the set of unfalsified by data coefficients of the transfer function of the nominal model remains the same a priori polytope.  Nevertheless, the adaptive control is able to guarantee, with the prescribed accuracy, the same optimal asymptotic upper bound on the plant output as for the plant with known parameters \cite{sok01b}. The optimal problem is considered within the $\ell_1$-theory of robust control associated with the $\ell_\infty$ signal space and bounded disturbances. Basic results on necessary and sufficient conditions for robust stability and robust performance in the $\ell_1$ setting were obtained for multidimensional LTI systems under zero initial conditions and structured norm bounded perturbations  that can be time varying or nonlinear \cite{kp91,kp93}. Necessary and sufficient conditions of robust stability and steady-state performance were proven for systems under fading or finite memory perturbations and fixed external signal \cite{kham95}. Then the representation for the worst-case norm of system output under bounded disturbances and structured perturbations was derived  \cite{kham97}. Since the models of fading and finite memory perturbations are not verifiable by data, the verifiable model of bounded memory perturbations was introduced for purposes of adaptive control and the representation for the worst-case steady-state  norm of the system output was obtained for general case of system with bounded disturbance, structured norm-bounded perturbations and additional fixed input (tracking signal) \cite{sok01a,sok99}.

The main ideas for synthesis of adaptive optimal control are, first, in the use of set estimates of the unfalsified by data parameters, the norms of the disturbance and perturbations including,  and, second, in the computation of current optimal vector estimates with the use of the control criterion as the identification criterion. Note that these ideas can not be used  for adaptive optimal robust control in the $H_\infty$ setting. Indeed, only auxiliary problem of assessing the quality of a given, or estimated, model, which was called  in 1997 ``A central issue in system identification''  \cite{ljungguo97},  remains an open issue to present time \cite{lam17}. The adaptive optimal controllers in the $\ell_1$ setting, based on these ideas, were proposed for systems with bounded disturbance \cite{sok85scl} and then for systems with additional perturbations \cite{sok96,sok01b}. These adaptive optimal controllers assumed very high computer power for computation of current optimal estimates, but  demonstrated the theoretical solvability of the problem of \textit{adaptive optimal robust control} of non-identifiable systems and, thus, showed the maximum capability of feedback in the $\ell_1$ setting. It was noted recently \cite{guo20}, ``In fact, the understanding of the maximum capability of feedback can encourage us in improving the controller design to reach or approach the maximum capability, and may help us in alleviating the workload of modeling and identification''.

The purpose of the present paper is to present a computationally tractable solution of the optimal problem described above. For computation of the current optimal estimates it is necessary to minimize the nonconvex control criterion  over the current set estimate of unfalsified by data parameters. This set estimates are described by data-based linear inequalities and a nonconvex inequality that describes the set of robustly stabilizable systems. We introduce a new estimated parameter instead of the two estimated norms of coprime factor perturbations. This change of the estimated parameters, under some nonrestrictive additional assumption about the total disturbance in the controlled system, allows to describe this system as a system with the bounded disturbance and the norm bounded output uncertainty. For such a system, the control criterion becomes a linear-fractional function of the norms of the additive disturbance and the output uncertainty, and the set of parameters of stabilizable systems is described by a linear inequality \cite{ait21}.  Since linear-fractional programming is reduced to linear programming \cite{boyd04}, the computation of the optimal current estimates becomes a computationally tractable problem. The number of linear inequalities in the description of the set estimates of unknown parameters is ensured to be bounded due to the use of a dead zone under their updating. The computational tractability of the proposed adaptive optimal robust control is illustrated by simulations for the system with 10 unknown parameters. Simulations with the RLS estimation algorithm, which has no proven results on the stability of closed loop system with this algorithm, are presented to illustrate its unsuitability to the adaptive optimal control in the $\ell_1$ setting.

The main contributions of this paper are as follows.
\begin{enumerate}
 \item 
The problem of adaptive robust stabilization of a discrete-time minimum-phase plant under time-varying or nonlinear coprime factor perturbations and bounded external disturbance is considered and solved in the $\ell_1$-optimal setting. The coefficients of the transfer function, the norms of coprime factor perturbations, and the norm of the additive external disturbance are assumed to be unknown. The control criterion is the worst-case upper upper bound on the plant output in steady-state. The use of the dead zone under estimates updating ensures the convergence of set and vector estimates in a finite time.

\item
Despite the non-identifiability of the unknown parameters,  the proposed adaptive controller guarantees, with the prescribed accuracy, the same optimal parameter dependent upper bound for the output of adaptive system as the optimal controller for the plant with known parameters. The accuracy of solving the optimal problem is determined by the choice of the size of the dead zone under estimates updating and can be regulated on-line depending on the current optimal estimates of the unknown parameters. 

\item
The norms of admissible coprime factor perturbations are not the same for different plants and can be arbitrary from any compact subset of the open set of the norms, for which the unknown plant is robustly stabilizable. Thus, the presented adaptive control realize the maximum capability of feedback with respect to both the control criterion  and the region of the admissible norms of coprime factor perturbations.

\item
In contrast to traditional adaptive systems, the current optimal and finite-time convergent estimates are verified (validated) by data on-line and are consistent with data in steady-state. The values of the control criterion for the current optimal estimates of unknown parameters give current optimal, unfalsified by data, upper bounds of the system output. The steady-state unfalsified value of the control criterion can be considerably less than the optimal value of the control criterion for the controlled plant depending on specific bounded disturbance and perturbations. Moreover, the unfalsified values of the control criterion can be considered as a criterion for verification of a priori assumptions about the controlled system.

\end{enumerate}

This paper is organized as follows. Problem statement is given in section \ref{ps}. Non-identifiability of estimated parameters is explained in section \ref{secnonid}. Motivation (or rather informal explanation of the necessity) of the set-member estimation and the use of the control criterion as the identification criterion is given in section \ref{secoe}. The main results on the convergence of estimates in a finite-time and the optimality of adaptive stabilization are presented in section \ref{secaos}. A piori assumptions, problems of model verification, and on-line choice of the dead-zone parameter under estimates updating are the topis of section \ref{secmv}. The computational tractability of the proposed adaptive control is illustrated in section \ref{secsim} by simulations for a plant with 10 unknown parameters.  Section \ref{con} concludes the paper.

\textit{Notation:}

\noindent
$|\varphi|$ -- the euclidean norm of vector $\varphi\in\mathbb{R}^n$. \\
$\ell_e$ -- space of real sequences $x=(\cdots,x_{-1},x_0,x_1,\cdots)$. \\
$x_s^t=(x_s,x_{s+1},\ldots,x_t)$ for $x\in\ell_e$.\\
$|x_s^t|=\max_{s\le k \le t} |x_k|$ for $x\in\ell_e$. \\
$\ell_\infty$ -- normed space of bounded real sequences, $\|x\|_{\ell_\infty}=\sup_t |x_t|$ for $x\in\ell_\infty$.\\
$\ell_1$ -- normed space of absolutely summable real sequences, $\|x\|=\sum_{k=0}^{+\infty}|x_t|$ for $x\in\ell_1$.\\
$\|x\|_{ss}=\limsup_{t\to+\infty} \ |x_t|$ for $x\in \ell_e$.\\
$\|G\|=\sum_{k=0}^{+\infty}|g_k|=\|g\|_{\ell_1}$  -- the induced norm of a stable causal linear time invariant system $G:\ell_\infty\to\ell_\infty$ associated with the transfer function $G(\lambda)=\sum_{k=0}^{+\infty} g_k\lambda^k$.\\


\section{Problem statement}
\label{ps}
\subsection{System description}

Consider a discrete-time single-input single-output (SISO) system described by
\begin{equation}
\label{plant}
a(q^{-1})y_{t+1}=b(q^{-1})u_t+v_{t+1} \quad t=0,1,2,\ldots , 
\end{equation}
\[
 a(q^{-1})=1+a_1 q^{-1}+\ldots+a_n q^{-n}\,,
\]
\[
 b(q^{-1})=b_1+b_2 q^{-1}+\ldots+b_m q^{1-m}\,,
\]
where $q^{-1}$  is the backward shift operator  ($q^{-1}y_t=y_{t-1}$), $y_t \in \mathbb R$, $u_t \in \mathbb R$, and $v_t \in \mathbb R$ are, respectively, the measured output, control input, and total disturbance in the system.  The initial values $y_{1-n},\cdots,y_0$ are arbitrary. . We set $y_k=0$ for all $k\le -n$ and $u_k=0$ for all $k<0$ to simplify the subsequent presentation.  A priori information about the system is as follows.

\textbf{Assumption A1.} The coefficients of the polynomials $a$ and $b$, which describe the nominal model, are in a (bounded) polytope $\Xi$, 
\[
\xi:=(a_1,\ldots,a_n,b_1,\ldots,b_m)^T
\in\Xi=\{ \hat\xi \ | \ P\hat\xi\ge p  \} \subset \mathbb{R}^{n+m}
\]
where the matrix $P\in\mathbb{R}^{l\times(n+m)}$ and the vector $p\in\mathbb{R}^l$ are known. It is assumed that
$b_1 \ne 0$  and the roots of $b(\lambda)$ are outside of the unit disk $\{z\in\mathbb{C} \ | \ |z|\le 1\}$ for any $\xi\in\Xi$, that is, the nominal models are minimum phase for any $\xi\in\Xi$ .

\textbf{Assumption A2.} The total disturbance $v$ is described by
\begin{align}
\label{v1}
& v_t =\delta^w w_t+ \delta^y\Delta^1(y)_t+\delta^u \Delta^2(u)_t\,, \\
\label{v2}
&  \|w\|_{\ell_\infty}\le 1, \quad
 |\Delta^1(y)_t| \le |y_{t-\mu}^{t-1}|, \quad 
 |\Delta^2(u)_t|\le |u_{t-\mu}^{t-1}|\,.
\end{align}
In (\ref{v1}), $w\in\ell_\infty$ is the normalized external disturbance,  $\delta^w$ is the upper bound of the external disturbance $\delta^w w$, the nonnegative $\delta^y$ and $\delta^u$ are, respectively, the gains (the induced $\ell_\infty$ norms) of output  and control perturbations $\delta^y\Delta^1$ and $\delta^u\Delta^2$. These perturbations are also called coprime factor perturbations of the system transfer function. The operators  $\Delta^1:\ell_\infty\to\ell_\infty$ and $\Delta^2:\ell_\infty\to\ell_\infty$ are normalized linear time-varying or nonlinear strictly causal operators with a bounded memory $\mu$ \cite{sok01a}.

\textbf{Assumption A3.} The parameter vector $\theta=(\xi^T,\delta^w,\delta^y,\delta^u)$  of the system (\ref{plant})-(\ref{v2}) is unknown.

 Another nonrestrictive technical a priori assumption on robust stabilizability of controlled system will be formulated at the end of subsection \ref{pf} before the strict formulation of the problem.
 
The problem under consideration is to design adaptive control that provides as small as possible upper bound for  the control criterion
\begin{equation}
\label{cc}
 \sup_{v\in V} \ \limsup_{t\to+\infty} \ |y_t|\,,
\end{equation}
where $V$ is the set of total disturbances $v$ satisfying Assumption \textbf{A2}. Strict formulation of the problem is given in the subsection \ref{pf}.


\subsection{Robust perfomance of optimal closed loop system with known parameters}

Consider the system in (\ref{plant}) with the known vector of coefficients $\xi$. The controller of the form
\begin{equation}
\label{optcont}
 b(q^{-1})u_t=(a(q^{-1})-1)y_{t+1}
\end{equation}
ensures the equality
\begin{equation}
\label{yopt}
 y_{t+1}=v_{t+1}=\delta^w w_{t+1}+ \delta^y\Delta^1(y)_{t+1}+\delta^u \Delta^2(u)_{t+1} \quad \forall \ t
\end{equation}
for the outputs $y_{t+1}$ of the system in (1) and, therefore, is optimal for the control criterion (\ref{cc}) in view of the unpredictability of $v_{t+1}$. Introduce notation
\[
 G^\xi(\lambda)=\frac{(a(\lambda)-1)\lambda}{b(\lambda)}=\sum_{k=0}^{+\infty}\, g^\xi_k\,\lambda^k
\]
for stable transfer function of the optimal controller in (\ref{optcont}) where $g^\xi$ is the impulse response of $G^\xi$ and
\[
 \|G^\xi\|=\sum_{k=0}^{+\infty}\, |g^\xi_k|=\|g^\xi\|_{\ell_1}\,,
\]
 Define the control criterion in the form
\begin{equation}
 \label{Jmu}
 J_\mu(\theta)=\sup_{v\in V} \ \limsup_{t\to+\infty} \ |y_t|
\end{equation}
where $y$ is the output of the optimal closed loop system (\ref{plant}) and (\ref{optcont}) and $V$ is the set of total disturbances $v$ satisfying Assumption \textbf{A2}.
 
The closed-loop system (\ref{plant}) and (\ref{optcont}) is called \emph{robustly stable} if $ J_\mu(\theta)<+\infty$. Robust performance of this system is described by the following theorem.

\begin{theorem}
 \label{optstab}
Following statements hold for the optimal closed loop system  (\ref{plant}) and (\ref{optcont}). 

1. The optimal closed loop system  (\ref{plant}) and (\ref{optcont})  with the perturbations memory $\mu=+\infty$ is robustly stable if and only if
\begin{equation}
 \label{rs}
 \delta^y+\delta^u\|G^\xi\|<1\,.
\end{equation}
and 
\begin{equation}
 J(\theta):=J_{+\infty}(\theta)=\frac{\delta^w}{1-\delta^y-\delta^u\|G^\xi\|}
\end{equation}
for the system with the zero initial data $y_{1-n}^0$.

2. For the system  (\ref{plant}) and (\ref{optcont}) with bounded memory perturbations ($\mu<+\infty$) and arbitrary initial data $y_{1-n}^0$,
\begin{equation}
 \label{Jmu1}
 J_\mu(\theta) \nearrow J(\theta) \quad (\mu\rightarrow +\infty)\,,
\end{equation}
where the sign $\nearrow$ denotes the monotone convergence from below.
\end{theorem}
\begin{proof}
The first statement of Theorem \ref{optstab} follows from Proposition 3\cite{kham97} applied to the system  (\ref{plant}) and (\ref{optcont}) while the second statement of Theorem \ref{optstab} follows from Theorems 6 and 7\cite{sok01a} applied to this system.
\end{proof}


 \subsection{Problem formulation}
 \label{pf}

 Before precise formulation of the optimal problem under consideration, we have to make some comments associated with the model of bounded memory perturbations in (\ref{v2}). Suchlike models of perturbations were used in problems of robust adaptive control since the late 1980s in more conservative forms (e.g. $|v_t|\le \delta^w+\delta\max{|y_{t-\mu}^{t-1}|, |u_{t-\mu}^{t-1}|}$ with the known $\delta^w, \delta$\cite{wey94}). Basic results on robust stability and robust performance of systems in the $\ell_1$ setting were related to  systems with infinite memory perturbations ($\mu=+\infty$) and only zero initial data\cite{kp91,kp93}, and, therefore, can not be applied for identification and adaptive control. These results were extended to steady-state  performance of systems with arbitrary initial data by assuming finite/fading memory perturbations instead of infinite memory perturbations\cite{kham95,kham97}. However, the finite/fading memory perturbations are not verifiable by measurement data because one can not test on-line whether a sequence of real numbers is finite or converges to zero. In contrast to this, the model of bounded memory perturbations is verifiable by data and, in addition, makes possible data based estimation of their gains (norms). Since, in view of (\ref{Jmu1}), the value of $J(\theta)$ is a tight upper bound of the worst-case value of $\|y\|_{ss}$ for large values of  $\mu$, we will treat $J(\theta)$  as the control criterion in the considered problem. 

 In order to formulate strict results on adaptive optimal control, another technical a priori assumption is used.
  
\textbf{Assumption A4.} 
 The unknown vector $\theta$ of the system parameters satisfies the inequality
 \begin{equation}
 \label{assA4}
   \delta^y+\delta^u\|G^\xi\|\le \bar\delta<1
 \end{equation}
with a known $\bar\delta$.

In fact, the value of $\bar\delta$ in (\ref{assA4}) is chosen by the controller designer and can be taken arbitrarily close to 1. Assumption \textbf{A4} is nonrestrictive. If $\delta^w\ne 0$ and $\delta^y+\delta^u\|G^\xi\|$ is very close to 1, then $J(\theta)\to+\infty$ as $\bar\delta\to 1$ and such models are useless in practice.

\textbf{Problem formulation.} We are interested in the synthesis of feedback of the form $u_t=U_t(y_{1-n}^t,u_0^t)$ that provides, with the prescribed accuracy, the inequality
\begin{equation}
 \label{problem}
 \|y\|_{ss}= \ \limsup_{t\to+\infty} \ |y_t| \le J(\theta) \,.
\end{equation}


\section{Nonidentifiability}
\label{secnonid}
 
 In this subsection,  two simple statements illustrate the complexity of the stated problem in view of the consistency of any $\hat\xi\in\Xi$ with measurement data on any finite time interval under a priori Assumptions \textbf{A1}-\textbf{A4}. Introduce notation
 \[
  p^y_t=|y_{t-\mu}^{t-1}|, \quad p^u_t=|u_{t-\mu}^{t-1}|\,.
 \]
\begin{proposition}
\label{proproi}
If some estimate $\hat\theta=(\hat\xi^T,\hat\delta^w,\hat\delta^y,\hat\delta^u)^T$, $ \hat\xi\in\Xi$, $\hat\delta^w\ge 0$, $\hat\delta^y\ge 0$, $\hat\delta^u\ge 0$, satisfies the inequalities
\begin{equation}
 \label{roi}
 |\hat a(q^{-1})y_{t+1}-\, \hat b(q^{-1})u_t|\le \hat\delta^w+\hat\delta^y p^y_{t+1}+\hat\delta^u p^u_{t+1}
\end{equation}
for all $t\ge 0$ then the plant (\ref{plant}) with the parameter vector $\hat\theta$ satisfies the equation (\ref{plant}) and a priori Assumptions \textbf{A1,A2} for all $t\ge0$.
\end{proposition}
\begin{proof}
Define $\hat v_{t+1}=\hat a(q^{-1})y_{t+1}-\hat b(q^{-1})u_t$. Then the plant (\ref{plant}) associated with the parameter vector $\hat\theta$ and the total disturbance $\hat v$ satisfy the inequality (\ref{plant})with the total disturbance $\hat v$ and 
\begin{equation}
 \label{v3}
|\hat v_{t+1}|\le \hat\delta^w+\hat\delta^y p^y_{t+1}+\hat\delta^u p^u_{t+1}\,.
\end{equation}
According to Lemma 1\cite{sok01a},  for any sequences $x$ and $z$ from $\ell_e$ the inequality $|x_t|\le |z_{t-\mu}^{t-1}|$ is equivalent to the existence of linear time varying  strictly causal operator  operator $\Delta:\ell_e\rightarrow\ell_e$ with a bounded memory $\mu$ such that $x=\Delta z$. This statement and the inequality (\ref{v3}) imply that the total disturbance $\hat v$ in the plant with the parameter vector $\hat\theta$ can be presented in the form (\ref{v1}) and this plant satisfies  Assumptions \textbf{A1,A2}.
\end{proof}
 
 It follows from Proposition \ref{proproi}, that full information about the unknown vector $\theta$ at any time $t$ and for any control inputs $u_0^{t-1}$ is of the form
 \begin{align*}
 \theta\in \Theta_t=\{ \ \hat\theta\in \Theta_0 \ \big | \ &  |\hat a(q^{-1})y_{k+1}-\hat b(q^{-1})u_k|\le \\
 & \hat\delta^w+\hat\delta^y p^y_{k+1}+\hat\delta^u p^u_{k+1} \ \forall k<t \ \}\,,
\end{align*}
where
\[
 \Theta_0=\{ \ \hat\theta \ \big |
 \hat\xi\in\Xi\,, \hat\delta^w\ge 0\,,  \hat\delta^y\ge 0\,,  \hat\delta^u\ge 0\,, \hat\delta^y+\hat\delta^u\|G^{\hat\xi}\|\le \bar\delta \ \}\!\,.
\]
The set $\Theta_0$ is the prior set of feasible values of $\theta$ and $\Theta_t$ is the set of estimates $\hat\theta$ unfalsified by (or, equivalently, compatible with) the measurement data $y_0^t$, $u_0^{t-1}$ and a priori Assumptions  \textbf{A1-A4}.
 
\begin{proposition}[Non identifiability of $\xi$]
 \label{nonid}
 For any control inputs $u_0^{t-1}$, any $\hat\xi\in\Xi$, and any nonnegative $\hat\delta^y,\hat\delta^u$,
\[
 \hat\theta=(\hat\xi^T,\hat\delta^w,\hat\delta^y,\hat\delta^u)^T\in \Theta_t
\]
for all sufficiently large $\hat\delta^w$.
\end{proposition}
\begin{proof}
 Proposition \ref{nonid} clearly follows from the monotone convergence of the right-hand side of the inequalities (\ref{roi}) to $+\infty$ as $\delta^w\to+\infty$.
\end{proof}

Proposition \ref{nonid} means that the set of compatible with data vectors $\hat\xi$ is not reduced with getting new data and remains the same a priori set $\Xi$.

\section{Optimal estimation} 
\label{secoe}
 
In 1960th, Prof. V.A. Yakubovich proposed the method of recurrent objective inequalities  for synthesis of adaptive control of dynamical systems under bounded disturbance with the known upper bound $\bar\delta^w$. The idea of the method was to compute an estimate $\hat\xi$ that would satisfy the inequalities 
$
 |\hat a(q^{-1})y_{t+1}-\, \hat b(q^{-1})u_t|\le \bar\delta^w
$
for all sufficiently large $t$. Various algorithms providing convergence of estimates of $\xi$ in a finite time have been proposed and applied to problems of adaptive stabilization \cite{bonyak92}. Note at first that identifying the vector $\xi$ alone, without evaluating $\delta^w, \delta^y, \delta^u$, is not sufficient to solve the stated optimal problem.  Further, if some extended estimate $\hat\theta$ satisfies inequalities (\ref{roi}), which can play the role of the objective inequalities in problem (\ref{problem}), for all sufficiently large $t$ and we apply the optimal controller corresponding to this estimate, then Theorem \ref{optstab} and Proposition \ref{proproi} guarantee the inequality
\[
 \|y\|_{ss}\le J(\hat\theta)\,.
\]
This inequality is insufficient for solution of the stated problem (\ref{problem}) and we need the additional inequality
\begin{equation}
 \label{Jineq}
 J(\hat\theta)\le J(\theta)\,.
\end{equation}
Since any vector in the set $\Theta_t$ is compatible with data and a priori information and, therefore, can be the unknown ``true'' vector $\theta$ of the plant (\ref{plant}), the inequality (\ref{Jineq}) dictates the choice of the control criterion $J$ as the identification criterion:
\begin{equation}
 \label{optest}
 \theta_t=\argmin_{\hat\theta\in \Theta_t} \ J(\hat\theta)=\argmin_{\hat\theta\in \Theta_t} \
 \frac{\hat\delta^w}{1-\hat\delta^y-\hat\delta^u\|G^{\hat\xi}\|}\,.
\end{equation}
On-line solution of the optimal estimation problem (\ref{optest}) with the prescribed accuracy is difficult because the control criterion $J$ and the constraint (\ref{assA4}) in Assumption \textbf{A4} are nonconvex. Moreover, even the approximate computation of $\|G^{\hat\xi}\|$ alone needs special computations\cite{pic13,san98} in view of no analytical representation of $\|G^{\hat\xi}\|$.
We will transform the problem (\ref{optest})  into a linear fractional problem that can be solved on-line.  This is achieved by a special change of the vector of unknown parameters that reduce, under additional assumption about the total disturbance $v$, the plant model (\ref{plant}) to a model with output only perturbation.

 The control sequence $u$ in the optimal closed loop system (\ref{plant}) and (\ref{optcont}) satisfies the inequalities
 \begin{align}
 \label{uub0}
 |u_t|=\left|G^\xi(q^{-1})y_t\right|=\left|\sum_{k=0}^{t+n-1}\, g^\xi_k\,y_{t-k}\right|
 \le |\sum_{k=0}^{\bar\mu}\, g^\xi_k\,y_{t-k}|+ \\ \nonumber
 |\sum_{k=\bar\mu+1}^{t+n-1}\, g^\xi_k\,y_{t-k}|
 \le  \sum_{k=0}^{\bar\mu}\, |g^\xi_k|\,|y_{t-k}|+\sum_{k=\bar\mu+1}^{t+n-1}\, |g^\xi_k|\,|y_{t-k}|\,.
\end{align}
Taking into account that $\sum_{k=0}^{\bar\mu}\, |g^\xi_k| \le \sum_{k=0}^{+\infty}\, |g^\xi_k|=\|G^\xi\|$ and $\sum_{k=\bar\mu+1}^{+\infty}\, |g^\xi_k|\rightarrow 0$ as $\bar\mu\rightarrow+\infty$, we shall assume that for some sufficiently large natural number $\bar\mu>\mu$, chosen by the designer, and any $\xi\in\Xi$ it holds 
 \begin{equation}
\label{uub}
 |u_t|=| \sum_{k=0}^{t+n-1}\, g^\xi_k\,y_{t-k}| \le \|G^\xi\| |y_{t+\mu-\bar\mu}^{t}|\,.
\end{equation}
Since $y=v$ in the optimal system  (\ref{plant}) and (\ref{optcont}), the assumption (\ref{uub}) is in fact the assumption 
\begin{equation}
\label{vass}
 | \sum_{k=0}^{t+n-1}\, g^\xi_k\,v_{t-k}| \le \|G^\xi\| |v_{t+\mu-\bar\mu}^{t}| \quad \forall \xi\in \Xi \,.
\end{equation}
about the total disturbance $v$. This assumption excludes total disturbances $v$ that maximize, from time to time, the absolute value of the control input $u_t$ in the optimal system. The assumption (\ref{vass}) actually becomes less and less restrictive with increasing $\hat\mu$ because the set of the total disturbances $v$ that don't satisfy the inequalities (\ref{vass}) converges to the empty set as $\bar\mu\to +\infty$. One can consider (\ref{vass}) as an additional assumption that the total disturbance $v$ is unintentional or casual. 

It follows now from (\ref{uub})
\begin{equation}
\label{puub}
 p^u_t=\max_{t-\mu\le k < t}|u_k|\le \|G^\xi\| |y_{t-\bar\mu}^{t-1}|\,.
\end{equation}
Introduce notation 
\begin{equation}
\label{pt}
 p_{t+1}=|y_{t+1-\bar\mu}^{t}| \,.
\end{equation}
It follows now from the plant equation (\ref{plant}), a priori Assumptions, and (\ref{puub})
\begin{align}
\label{roi1}
 | a(q^{-1})y_{t+1}- b(q^{-1})u_t| & \le \delta^w+\delta^y p^y_{t+1}+ \delta^u p^u_{t+1} \\ 
 & \le  \delta^w+(\delta^y+\delta^u \|G^{\xi}\|)p_{t+1}\,.\nonumber
\end{align}
Introduce new parameters to be estimated
\begin{equation}
\label{zeta}
\delta=\delta^y+\delta^u\|G^\xi\|\,, \quad  \zeta= (\xi,\delta^w,\delta)^T\,.
 \end{equation}
In this notation, the  inequalities (\ref{roi1}) take the form
\begin{equation}
 \label{roinew}
 | a(q^{-1})y_{t+1}- b(q^{-1})u_t|\le \delta^w+\delta p_{t+1} \ \ \forall\, t 
\end{equation}
and
\[
 I(\zeta):=\frac{\delta^w}{1-\delta}=J(\theta)\,.
\]
We can consider the inequalities (\ref{roinew}) as the inequalities (\ref{roi}) for the plant with the parameter vector
$\hat\theta=(\xi^T,\delta^w,\delta,0)^T$. In view of the Proposition \ref{proproi}, the inequalities  (\ref{roinew}) imply that the sequence $y$ can be considered as the output of the plant (\ref{plant}) with the parameter vector
$\hat\theta$ and then
\begin{equation}
 \label{I=J}
 J(\hat\theta)=\frac{\delta^w}{1-\delta}=I(\zeta)
\end{equation}
by Theorem \ref{optstab}.
After the described change of parameters to be estimated, the problem of optimal estimation (\ref{optest}) becomes the problem of linear fraction programming
\begin{equation}
 \label{optestnew}
 \zeta_t=\argmin_{\hat\zeta= (\hat\xi,\hat\delta^w,\hat\delta)^T\in S_t}  J(\hat\zeta)=\argmin_{\hat\zeta\in S_t} \
 \frac{\hat\delta^w}{1-\hat\delta}\,,
\end{equation}
where
\[
 S_t=\{  \hat\zeta\in Z_0 \ \big | 
 |\hat a(q^{-1})y_{k}-\hat b(q^{-1})u_k|\le \hat\delta^w+\hat\delta p^y_{k} \ \forall k\le t+1  \}\,,
\]
\begin{equation}
\label{Z_0}
 Z_0=\{ \ \hat\zeta \ =(\hat\xi^T,\hat\delta^w,\hat\delta)^T \ \big | \
 \hat\xi\in\Xi\,, \ \hat\delta^w\ge 0\,,  \ 0\le\hat\delta \le \bar\delta \ \}\,.
\end{equation}
Thus, the nonconvex control criterion $J(\theta)$ takes the form of the linear fractional criterion $I(\zeta)$, nonconvex a priori constraint (\ref{assA4}) in Assumption \textbf{A4} is transformed into the linear constraint $\hat\delta \le \bar\delta$ and information inequalities (\ref{roi}) are replaced by the inequalities (\ref{roinew}). As a result, the problem of optimal estimation (\ref{optest}) is transformed to the linear fractional problem  (\ref{optestnew}), which is reducible in a standard way to a linear programming\cite{boyd04}.


 \section{Adaptive optimal stabilization}
 \label{secaos}

The number of inequalities in the description of the sets $S_t$ of compatible with data vectors $\hat\zeta$ can grows without limit with time. To avoid this, we will use upper set estimates $Z_t\supset S_t$ and  a dead zone under updating set estimates that guarantees the convergence of the estimates in finite time. Choose a small real number $\varepsilon>0$, the parameter of the dead zone. The less is $\varepsilon>0$, the closer to (\ref{problem}) will be the adaptive control performance. 

We define the estimation algorithm and the adaptive controller as follows. Choose a natural number $\bar\mu\ge 2\mu$.  At every time instant $t$, it is computed a polyhedral estimate $Z_t$ and a vector estimate $\zeta_t$
\[
 \zeta_t=(\xi^T_t,\delta^w_t,\delta_t)^T
\]
of the unknown vector $\zeta=(\xi,\delta^w,\delta)^T$. Define the initial estimate $Z_0$ by (\ref{Z_0}) and the initial $\zeta_0=(\xi_0^T,0,0)^T$ with arbitrary $\xi_0\in\Xi$. The control input $u_t$ at the time instant $t$ is computed in two steps. At first, a preliminary value of $u_t$ is computed by the adaptive controller
\begin{equation}
 \label{ac}
 b^t(q^{-1})u_t=(a^t(q^{-1})-1)y_{t+1} \,.
\end{equation}
At the second step, $u_t$ is corrected, if necessary, as follows
\begin{equation}
 \label{cut}
 u_t:=\sign(u_t)\|G^{\xi_t}\| |y_{t+\mu-\bar\mu}^t|\,, \mbox{ \ if \ } |u_t| > \|G^{\xi_t}\| |y_{t+\mu-\bar\mu}^t|\,.
\end{equation}
The cutting (\ref{cut}) guarantees that the inequalities
\begin{equation}
 \label{uub1}
 |u_t| \le \|G^{\xi_t}\| |y_{t+\mu-\bar\mu}^t|
\end{equation}
hold for any $t\ge 0$. 

The polyhedral estimates $Z_t$ and the vector estimates $\zeta_t$ are updated as follows. Having measured the value of $y_{t+1}$, define
\begin{equation*}
 \begin{array}{lll}
 &\varphi_t:=(-y_t,-y_{t-1},\ldots,-y_{t-n+1},u_t,\ldots,u_{t-m+1})^{T}\,, \\  
 &\eta_{t+1}:=\sign\, (y_{t+1}-\varphi_t^{T}\xi_t)\,, \\
 &\psi_{t+1}:=(\eta_{t+1}\varphi_t^{T},1,p_{t+1})^{T}\,, \quad
 \nu_{t+1}:=\eta_{t+1}y_{t+1}\,.
 \end{array}
\end{equation*}
In this notation, the inequality (\ref{roinew}) with respect to $\zeta_t$ is equivalent to the inequality
\begin{equation}
\label{sokroi2}
 \psi_{t+1}^{T}\zeta_t \ge\nu_{t+1}\,.
\end{equation}
Define
\begin{equation}
\label{zetat1}
Z_{t+1}:=Z_t\,, \quad
\zeta_{t+1}:=\zeta_t, \mbox{ \ if \ } \psi_{t+1}^{T}\zeta_t\ge \nu_{t+1}-\varepsilon|\psi_{t+1}|\,.
\end{equation}
Otherwise
\begin{equation}
\label{zetat2}
Z_{t+1}:=Z_t \cap \Omega_{t+1}\,, \quad
\Omega_{t+1}:=\{ \ \hat\zeta \ \big | \ \psi_{t+1}^{T}\hat\zeta \ge\nu_{t+1} \}\,,
\end{equation}
\begin{equation}
\label{zetat3}
\zeta_{t+1}:=\argmin_{\hat\zeta\in Z_{t+1}} \ I(\hat\zeta)\,.
\end{equation}

The described estimation algorithm has a simple geometric interpretation. The estimate $Z_t$  is updated according to (\ref{zetat2}) if and only if the distance from the vector $\zeta_t$ to the halfspace $\Omega_{t+1}$ is greater than the dead zone parameter $\varepsilon$. In this case, the polyhedral estimate $Z_t$ is supplemented with the inequality $\psi_{t+1}^{T}\hat\zeta \ge\nu_{t+1}$ that defines the halfspace $\Omega_{t+1}$ in (\ref{zetat2}).

 \begin{theorem}
 \label{aoc}
 Let the plant(\ref{plant}) with the unknown parameter vector $\theta=(\xi^T,\delta^w,\delta^y,\delta^u)^T$ satisfy Assumptions \textbf{A1-A4} and be controlled  by the controller (\ref{ac}), (\ref{cut}) with the estimation algorithm (\ref{zetat1})--(\ref{zetat3}) and the dead zone parameter $\varepsilon$ satisfying 
 \begin{equation}
  \label{epsub}
  0<\varepsilon < (1-\bar\delta)/(2+G_u)\,, \quad G_u=\sup_{\xi\in\Xi}\|G^\xi\|\,.
 \end{equation}
 
Then the following statements hold.

1) If the number of cuttings (\ref{cut}) in the closed loop system is finite and the parameters $\delta^y$ and $\delta^u$ satisfy inequality
\begin{equation}
 \label{de1}
 \delta^y+\delta^u G_u \le \bar\delta<1\,,
\end{equation}
then the set estimates $Z_t$ and the vector estimates $\zeta_t$ converge in finite time and
\begin{align}
\label{adstab}
  \limsup_{t\to+\infty} \ |y_t| \le I(\zeta_\infty^\varepsilon) <   I(\zeta_\infty)+K_{\zeta_\infty}\varepsilon & \le \\ \nonumber
  \frac{\delta^w}{1-\delta^y-\delta^u \max_t \|G^{\xi_t}\|}+K_{\zeta_\infty}\varepsilon & \le \\
  \frac{\delta^w}{1-\delta^y-\delta^u G_u}+K_{\zeta_\infty}\varepsilon & \,,  \nonumber
\end{align}
where $\zeta_\infty=(\xi_\infty^T,\delta^w_\infty,\delta_\infty)^T$ is the final value of $\zeta_t$, $\zeta_\infty^\varepsilon=(\, \xi_\infty^T,\, \delta^w_\infty+\varepsilon,\,\delta_\infty+\varepsilon(2+\|G^{\xi_\infty}\|))^T$, and
 \begin{equation}
  \label{Kinf}
  K_{\zeta_\infty}=\frac{1+\delta^w_\infty(2+\|G^{\xi_\infty}\|)}
  {(1-\delta_\infty-\varepsilon(2+\|G^{\xi_\infty}\|))^2}\,.
 \end{equation}
  
2) If the number of cuttings (\ref{cut}) in the closed loop system is finite and the control $u$ satisfies for all $t$ the inequalities
\begin{equation}
 \label{uoptbound1}
 |u_t| \le \|G^\xi\| |y_{t+1-\bar\mu}^t|\,,
\end{equation}
then the set estimates $Z_t$ and the vector estimates $\zeta_t$ converge in finite time and $I(\zeta_\infty)\le J(\theta)$ so that
\begin{equation}
\label{suboptac}
  \limsup_{t\to+\infty} \ |y_t| \le I(\zeta_\infty^\varepsilon) <
  J(\theta)+K_{\zeta_\infty}\varepsilon \,.
\end{equation}
\end{theorem}

\begin{proof}
To prove the first statement of Theorem \ref{aoc}, we show at first that the distance from $\zeta_t$ to the halfspace $\Omega_{t+1}$ is greater than $\varepsilon$ under each updating $\zeta_t$. In view of (\ref{zetat1}), $ \psi_{t+1}^{T}\zeta_t <  \nu_{t+1}-\varepsilon|\psi_{t+1}|$ under each updating. Taking into account that $\psi_{t+1}^{T}\hat\zeta \ge\nu_{t+1}$ for any $\hat\zeta\in\Omega_{t+1}$, we get the inequality 
\[
  \varepsilon|\psi_{t+1}| <  |\psi_{t+1}^T(\hat\zeta-\zeta_t)| \le |\psi_{t+1}| |\hat\zeta-\zeta_t| \quad \forall \hat\zeta\in\Omega_{t+1}
 \]
and  $|\hat\zeta-\zeta_t|>\varepsilon$ for any $\hat\zeta\in\Omega_{t+1}$. Then any updated estimate $Z_{t+1}$ has the empty intersection with the $\varepsilon$-neighborhood of the estimate $\zeta_t$. In view of the monotone decreasing of the polyhedral estimates $Z_t$, the $\varepsilon/2$-neighborhoods of the updated estimates $\zeta_t$ have the empty intersections for all $t$. Then the number of possible updates of the estimates $Z_t$ and $\zeta_t$ is finite if all estimates $\zeta_t$ are in a bounded set. Now we prove the boundedness of the estimates $\zeta_t$.

In view of Assumptions \textbf{A1, A2}, (\ref{cut}), and (\ref{pt}) we have for the closed loop adaptive system (\ref{plant}) and (\ref{ac})
\begin{align}
 \label{roi2}
 & | a(q^{-1})y_{t+1}-\,  b(q^{-1})u_t|\le \delta^w+\delta^y p^y_{t+1}+\delta^u p^u_{t+1}\le \\ \nonumber
 & \delta^w+\delta^y p^y_{t+1}+\delta^u (\max_t \|G^{\xi_t}\|) |y_{t+\mu-\bar\mu}^t| \le \\ \nonumber
 & \delta^w+(\delta^y+\delta^u G_u)p_{t+1}\,, 
\end{align}
where $G_u$ is defined in (\ref{epsub}). The second inequality is equivalent to the inequality  $\psi_{t+1}^{T}\bar\zeta \ge\nu_{t+1}$ with $\bar\zeta=(\xi^T,\delta^w,\delta^y+\delta^u \max_t \|G^{\xi_t}\|)^T$, that is $\bar\zeta\in\Omega_{t+1}$ and, consequently, $\bar\zeta\in Z_{t+1}$.  In vew of (\ref{zetat3}) we get the inequalities
\begin{equation}
\label{Jup}
 I(\zeta_t)\le I(\bar\zeta)=\frac{\delta^w}{1-\delta^y-\delta^u \max_t \|G^{\xi_t}\|}\le 
 \frac{\delta^w}{1-\delta^y-\delta^u G_u}
\end{equation}
 which implies the boundedness of $\zeta_t$ and, consequently, the convergence of the estimates $Z_t$ and $\zeta_t$ in a finite time $t_\infty$.
 
Since $\zeta_t=\zeta_\infty$ for all $t\ge t_\infty$, we have from (\ref{zetat1})
\[
 \psi_{t+1}^{T}\zeta_\infty\ge \nu_{t+1}-\varepsilon|\psi_{t+1}| \quad \forall t\ge t_\infty\,.
\]
These inequalities are equivalent to the inequalities
\[
 |a_\infty(q^{-1})y_{t+1}-b_\infty(q^{-1})u_t|\le \delta^w_\infty+\delta_\infty p_{t+1} +  \varepsilon|\psi_{t+1}|
\]
and then
 \begin{align}
   \label{oifin}
  & |a_\infty(q^{-1})y_{t+1}-b_\infty(q^{-1})u_t|\le 
  \delta^w_\infty+\delta_\infty p_{t+1}+ \\  \nonumber
  & \varepsilon(|y_{t-n+1}^{t}|+|u_{t-m+1}^{t}|+1+ p_{t+1})\le \\ \nonumber
  & \delta^w_\infty+\varepsilon+[\delta_\infty+\varepsilon(2+\|G^{\xi_\infty}\|)]p_{t+1}
\end{align}
 for all $t\ge t_\infty$. Due to Proposition \ref{proproi} and (\ref{oifin}) we can consider, for $t\ge t_{\infty}$, the output $y$ of the closed loop adaptive system as the output of the plant  of the form (\ref{plant}) with the parameter vector
\begin{equation}
 \label{thetateps}
 \theta_\infty^\varepsilon=(\, \xi_\infty^T,\, \delta^w_\infty+\varepsilon,\,
 \delta_\infty+\varepsilon(2+\|G^{\xi_\infty}\|), \, 0)^T\,.
\end{equation}
This vector corresponds to the final estimate $\zeta_\infty^\varepsilon=(\, \xi_\infty^T,\, \delta^w_\infty+\varepsilon,\,
 \delta_\infty+\varepsilon(2+\|G^{\xi_\infty}\|))^T$ of the estimation algorithm (\ref{zetat1})--(\ref{zetat3}) and $I(\zeta_\infty^\varepsilon)=J(\theta_\infty^\varepsilon)$. The inequalities (\ref{oifin}) and the inequality for $\varepsilon$ in (\ref{epsub}) imply the condition of robust stability (\ref{rs})
\begin{equation}
 \delta_\infty+\varepsilon(2+\|G^{\xi_\infty}\|)\le \bar\delta+\varepsilon(2+G_u)< 1
\end{equation}
\label{rsinf}
 for the plant with the parameter vector $\theta_\infty^\varepsilon$. Thus, under the condition of finite number of the cuttings in (\ref{cut}), the output $y$ can be considered, for all sufficiently large $t$, as the output of the optimal closed loop system of the form (\ref{plant}) and (\ref{optcont}) corresponding to the parameter vector $\theta_\infty^\varepsilon$ and in view of  Theorem \ref{optstab}
\[
 \limsup_{t\to+\infty} \ |y_t| \le I(\zeta_\infty^\varepsilon)=J(\theta_\infty^\varepsilon)=
 \frac{\delta^w_\infty+\varepsilon}
 {1-(\delta_\infty+\varepsilon(2+\|G^{\xi_\infty}\|))}.
\]
To estimate the difference $I(\zeta_\infty^\varepsilon)-I(\zeta_\infty)$, we use the inequality
\[
 \frac{C_1+\varepsilon_1}{C_2-\varepsilon_2}-\frac{C_1}{C_2}=
 \frac{C_2\varepsilon_1+C_1\varepsilon_2}{C_2(C_2-\varepsilon_2)} <
 \frac{\varepsilon_1+C_1\varepsilon_2}{(C_2-\varepsilon_2)^2}
\]
with the parameters
$
 C_1=\delta^w_\infty\,, \ 
 C_2=1-\delta_\infty\le 1\,, \ 
 \varepsilon_1=\varepsilon\,, \ 
 \varepsilon_2=\varepsilon(2+\|G^{\xi_\infty}\|)\,.
$
Then
\[
  I(\zeta_\infty^\varepsilon)-I(\zeta_\infty)< \frac{1+\delta^w_\infty(2+\|G^{\xi_\infty}\|)}
  {(1-\delta_\infty-\varepsilon(2+\|G^{\xi_\infty}\|))^2} \ \varepsilon=K_{\zeta_\infty}\varepsilon\,.
\]
This inequality is equivalent to the left inequality in (\ref{adstab}) with  $K_{\zeta_\infty}$ of the form (\ref{Kinf}). The right inequality in  (\ref{adstab}) followes now from (\ref{Jup}) and the convergence of $I(\zeta_t)$.

To prove the second statement of Theorem \ref{aoc}, we note that the condition in (\ref{uoptbound1}) implies similarly to (\ref{roi2})
\begin{align}
 \label{roi3}
 & | a(q^{-1})y_{t+1}-\,  b(q^{-1})u_t|\le \delta^w+\delta^y p^y_{t+1}+ \\ \nonumber
 & \delta^u \|G^{\xi}\| |y_{t+\mu-\bar\mu}^t| \le
 \delta^w+(\delta^y+\delta^u \|G^{\xi}\|)p_{t+1}\quad \forall t\,.
\end{align}
It means that the vector $\zeta=(\xi^T,\delta^w,\delta)^T$ with $\delta=\delta^y+\delta^u \|G^{\xi}\|$ satisfies the inequalities $\psi_{t+1}^T\zeta\ge\nu_{t+1}$ for all $t$. Then $\zeta\in \Omega_{t+1}$, $\zeta\in Z_{t+1}$  for all $t$ and $I(\zeta_t)\le  I(\zeta)$ in view of (\ref{zetat3}). As in the proof of the first statement, the inequality $I(\zeta_t)\le  I(\zeta)$ implies the convergence of estimates $Z_t$ and $\zeta_t$ in finite time. Then $I(\zeta_\infty)\le I(\zeta)$ and
\[
  \limsup_{t\to+\infty} \ |y_t| \le I(\zeta_\infty^\varepsilon)\,.
 \]
Taking into account the equality $I(\zeta)=J(\theta)$ (see (\ref{I=J})), we get  the second statement and Theorem \ref{aoc} is proven.
\end{proof}

\textit{Remark 1}.
Note that the convergence of the estimates $Z_t$ and $\zeta_t$ in finite time is proven in both statements of Theorem \ref{aoc} without the assumption of finite number of cuttings in the closed loop adaptive system. So this assumption is  exactly the assumption (\ref{vass}) with respect to the final estimates $\xi_\infty$, which can be any vector in a priori polyhedron $\Xi_0$ (indeed, $\xi_\infty=\xi_t=\xi$ for all $t$ if $\xi_0=\xi$). This assumption excludes ``deliberate'' total disturbances $v$ that almost maximize $|u_t|$ and the set of such disturbances comes to the empty set as $\bar\mu$ increases without limit.

\textit{Remark 2}. For the plant under consideration, the first statement of Theorem \ref{aoc} presents more advanced result achievable with the use of projection type estimation algorithm. Under additional assumption of known upper bound $\bar\delta^w$ on unknown norm $\delta^w$ of bounded disturbance, the projection algorithm described in \cite{sok01b} guarantees the inequality 
\begin{equation}
 \label{Jwey}
 \limsup_{t\to+\infty}\ |y_t| \le \frac{\bar\delta^w}{1-\bar\delta}+K_{\zeta_\infty}\varepsilon \,.
\end{equation}
This upper bound is the same for all admissible triples  $(\delta^w, \delta^y, \delta^u)$ and, therefor is considerable worse the upper bounds in (\ref{adstab}), which correspond to the ``true'' values of unknown parameters $\delta^w, \delta^y, \delta^u$ and decreases when any of these parameters decreases. Similar to (\ref{Jwey}) upper bound was presented in \cite{wey94} for autoregressive model with control delay and more conservative unstructured uncertainty.

\textit{Remark 3}. The second statement of Theorem  \ref{aoc} provides a solution to the problem (\ref{problem}) with the accuracy $K_{\zeta_\infty}\varepsilon$ under additional assumption (\ref{cut}). It must be emphasized that this assumption was satisfied in all  simulations performed with various disturbances $v$. The reason is in the fact that, first, the last inequality in (\ref{uub0}) is very rough and, second, optimal estimates (\ref{zetat3}) must have, if $I(\zeta_t)>0$, as small values of the component $\delta_{t+1}$ as possible to minimize the cost function $I(\hat\zeta)=\hat\delta^w/(1-\hat\delta)$. So it is difficult, if possible, to find a disturbance $v$ that can violate inequality in (\ref{cut}). At the same time, a strong mathematical proof of this inequality is an open problem. This assumption and the accuracy of the solution are commented and discussed  in the next section.

  
\section{Model verification in closed loop}
\label{secmv}

\subsection{Model verification under Assumption A4}
 The main and obvious advantage of the estimation algorithm (\ref{zetat1})-(\ref{zetat3}) is in the inequalities (\ref{suboptac}), which declare the approximate solution of the \textit{optimal} problem (\ref{problem}). Less obvious but equally important and unique benefits of the algorithm are that both the current estimates $\zeta_t$, the accuracy of the solution, and a priori assumptions are verified by data in closed loop. We will comment these benefits in more details. 
 
 The assumption of finite number of the cuttings (\ref{cut}) in the Theorem \ref{aoc} follows from the assumption (\ref{vass}) and present actually a weakened and verifiable by data version of the assumption (\ref{vass}), while the assumption (\ref{vass}) itself is not verifiable by data. Indeed, the control process for $t\ge t_\infty$ looks like the plant (\ref{plant}) with the parameter vector $\theta_\infty^{\varepsilon}$, defined in (\ref{thetateps}), is controlled by the optimal controller for this plant and the total disturbances $v^\infty_t= a_\infty(q^{-1})y_{t+1}- b_\infty(q^{-1})u_t$ in this plant cutisfies the a priori Assumptions 2 and 4.  Possible violations of the inequalities (\ref{puub}) at some time instants imply violations of the inequalities (\ref{vass}) for the plant at these time instants. Thus the inequalities (\ref{puub}) make  possible the verification of the assumption (\ref{vass}).  Note that he cuttings (\ref{cut}) were never observed in numerous simulations with various random and deterministic disturbances and perturbations.
 
The values of $\zeta_\infty$ and $K_{\zeta_\infty}$ are never known because any current estimate $\zeta_t$ can be falsified by future data. However, if there are no cuttings (\ref{cut}) from some time instant (or the cuttings occur less and less often), then, in view of the finite number of possible updates of the estimates $\zeta_t$, the intervals with the same estimate $\zeta_t$ become longer and longer. Therefore, the current estimates $\zeta_t$ that remain the same on long time intervals are  validated by data and the values of $I(\zeta_t)+K_{\zeta_t}\varepsilon$ become the best unfalsified and correct asymptotic upper bounds on the $|y_t|$.

It can happen that the value of $K_{\zeta_t}\varepsilon$ is greater than a desired accuracy of solution of the problem  (\ref{problem}). Consider this situation in detail. We have $I(\zeta_0)=0$ since $\delta^w_0=0$ and $I(\zeta_t)$ can remain zero on some time interval even if $I(\zeta)=J(\theta)\ne 0$. Note at first that the case $J(\theta)=0$, i.e. $\delta^w=0$, means no additive disturbance in the plant. The problem (\ref{problem}) is degenerate in this case because any controller, that robustly stabilizes the plant, is optimal with respect to the control criterion $J(\theta)$. If we know a priori that there is no external disturbance in the plant (\ref{plant}), then the control criterion (\ref{cc}) is of small interest. However, the controller (\ref{optcont}) clearly remains the best one for the plant (\ref{plant}) with respect to any reasonable control criterion because it provides the best possible dynamics, $y_t=v_t$ for all $t$, of the closed loop system. In this degenerate case, one could consider an optimal problem for another control criterion, e.g.
\[
 J(\theta)=\delta^y+\delta^u\|G^\xi\| \quad \Leftrightarrow \quad I(\zeta)=\delta\,.
\]
The optimal problem for this criterion can be solved in the adaptive setting (that is, for the plant with unknown parameters) with the use of the estimation algorithm  (\ref{zetat1})-(\ref{zetat3}) via simple eliminating the parameter $\delta^w$ from the vector $\zeta$.

Let us return to the problem (\ref{problem}) in the case of the plant with the external disturbance. In order to guarantee the desired absolute accuracy
\begin{equation}
 \label{E0}
 \limsup_{t\to+\infty} \ |y_t| \le  E
\end{equation}
for a chosen small positive $E$ while $I(\zeta_t)=0$, it suffices to guarantee the inequality 
\[
 I(\zeta_t^\varepsilon)= \frac{\varepsilon} {1-(\delta_t+\varepsilon(2+\|G^{\xi_t}\|))} \le E\, ,
\]
which is equivalent to the inequality
\begin{equation}
 \label{epsE}
 \varepsilon\le \varepsilon_t=\frac{1-\delta_t}{1+E(2+\|G^{\xi_t}\|)}E\,.
\end{equation}
Define the dead zone parameter $\varepsilon$ in (\ref{zetat1}) as $\varepsilon=\varepsilon_t$.  The convergence of the estimates $\zeta_t$ in finite time is preserved in view of the separation of $\varepsilon_t$ from zero
\begin{equation}
  \varepsilon_t \ge \frac{1-\bar\delta}{1+E(2+G_u)}E>0\,.
\end{equation}
 Then Theorem \ref{aoc} ensures the inequality (\ref{E0}) if $J(\zeta_t)=0$ for all $t$.
 
Consider now the nondegenerate case, when $I(\zeta_t)$ becomes nonzero. Let $I(\zeta_{t_*})>0$ and $I(\zeta_t)=0$ for $t<t_*$. Since the value of $J(\theta)$ is unknown a priori, it seems more natural to guarantee an approximate solution to the problem (\ref{problem}) in terms of the relative accuracy. Then the problem is to ensure the inequality
 \begin{equation}
  \label{relopt}
  \limsup_{t\to+\infty} \ |y_t| \le \kappa J(\theta)
 \end{equation}
 for a given $\kappa>1$.
 
For solution of the problem (\ref{relopt}), consider the following algorithm for updating $\varepsilon_t$. Choose any $\varkappa$, $\bar\delta<\varkappa<1$. The initial values of $\varepsilon_t$ for $t < t_*$ are the same as in (\ref{epsE}). For $t\ge t_*$ define
 \begin{equation}
  \label{epst}
  \varepsilon_{t}=\min\left(\frac{\varkappa-\bar\delta}{2+\|G^{\xi_t}\|}\,, \, \frac{(\kappa-1)I(\zeta_t)}{K_{\zeta_t}}\right)
\end{equation}

\begin{theorem}
 \label{aoc2}
  Let the plant(\ref{plant}) with the unknown parameter vector $\theta=(\xi^T,\delta^w,\delta^y,\delta^u)^T$ be controlled  by the controller (\ref{ac}), (\ref{cut}) and the dead zone parameter $\varepsilon$ in the estimation algorithm (\ref{zetat1})--(\ref{zetat3}) be defined by (\ref{epsE}) when $I(\zeta_t)=0$ and by (\ref{epst}) when $I(\zeta_t)\ne 0$. If the number of cuttings (\ref{cut}) in the closed loop system is finite, then the set estimates $Z_t$, the vector estimates $\zeta_t$, and the sequence \{$\varepsilon_t$\} converge in  finite time, the inequality (\ref{relopt}) holds if $I(\zeta_t)\ne 0$ for some $t$, and the inequality (\ref{E0}) holds if $I(\zeta_t)=0$ for all $t$.
\end{theorem}
\begin{proof}
 At first we prove that the number of the updates (\ref{zetat2}) is finite. The case $I(\zeta_t)=0$ for all $t$ was considered above. In the case $I(\zeta_{t_*})>0$ we have $I(\zeta_t)\ge I(\zeta_{t_*})$ for $t\ge t_*$ in view of increasing $I(\zeta_t)$.  In order to separate $\varepsilon$ from zero, note at first that
  \begin{equation}
 \label{sep1}
  \frac{\varkappa-\bar\delta}{2+\|G^{\xi_t}\|} \ge \frac{\varkappa-\bar\delta}{2+G_u}>0\,.
\end{equation}
The inequality $\varepsilon_t\le (\varkappa-\bar\delta)/(2+\|G^{\xi_t}\|)$, which follows from  (\ref{sep1}), implies
\begin{align}
 \label{rst}
 & \varepsilon_t\le \frac{\varkappa-\bar\delta}{2+\|G^{\xi_t}\|} \ \Rightarrow \
 \varepsilon_t(2+\|G^{\xi_t}\|) \le \varkappa-\bar\delta \ \Rightarrow \\ \nonumber
 & \delta_t+\varepsilon_t(2+\|G^{\xi_t}\|) \le \bar\delta+\varepsilon_t(2+\|G^{\xi_t}\|) \le \varkappa<1\,.
\end{align}
The last inequality in (\ref{rst}) guarantees the condition of robust stability (\ref{rsinf}) for the estimate $\theta_t^{\varepsilon_t}$ corresponding to the estimate $\zeta_t$.
It follows from (\ref{rst}), the representation (\ref{Kinf}) applied to $\zeta_t$, and the inequality $I(\zeta_t)\le I(\zeta)$ that
\begin{align}
 \label{Kmax}
 & K_{\zeta_t}\le \frac{1+\delta^w_t(2+\|G^{\xi_t}\|)}{(1-\delta_t-\varepsilon_t(2+\|G^{\xi_t}\|))^2} \le
 \frac{1+\delta^w_t(2+G_u)}{(1-\varkappa)^2} \le \\ \nonumber
 & \frac{\left(1+\frac{\delta^w_t(2+G_u)}{1-\delta_t}\right)}{(1-\varkappa)^2}\le 
 \frac{1+I(\zeta_t)(2+G_u)}{(1-\varkappa)^2}\le \\ \nonumber
 & \frac{1+I(\zeta)(2+G_u)}{(1-\varkappa)^2}=K_{max}.
\end{align}
Then we get for the right term in the minimization (\ref{epst})
\begin{equation}
 \label{sep2}
  \frac{(\kappa-1)I(\zeta_t)}{K_{\zeta_{t}}} \ge \frac{(\kappa-1)I(\zeta{t_*})}{K_{max}}>0\,.
\end{equation}
Now the separations (\ref{sep1}) and (\ref{sep2}) of $\varepsilon_t$ from zero imply, as in the proof of Theorem \ref{aoc}, the convergence of the estimates $\zeta_t$ and $Z_t$ in a finite time and, consequently, the convergence $\varepsilon_t\to\varepsilon_\infty$ in a finite time.
In view of no updates in (\ref{zetat2}), we have 
\[
  \psi_{t+1}^T\zeta_\infty \ge \nu_{t+1}-\varepsilon_\infty |\psi_{t+1}|
\]
for all sufficiently large $t$. Then it follows by  Theorem \ref{aoc}
\[
  \limsup_{t\to+\infty} \ |y_t| \le I(\zeta_\infty^{\varepsilon_\infty}) \le I(\zeta_{\infty})+K_{\zeta_\infty}\varepsilon_\infty \le 
  \kappa I(\zeta_{\infty})\,,
\]
where the last inequality follows from the inequality $\varepsilon_\infty\le (\kappa-1)I(\zeta_\infty)/K_{\zeta_\infty}$ provided by (\ref{epst}). Theorem \ref{aoc2} is proven.
\end{proof}

The choice of the dead zone parameter $\varepsilon$ in the Theorem \ref{aoc} needs the computation of $G_u$ in (\ref{epsub}) that can be a difficult problem. The described on-line computing the estimates $\varepsilon_t$ does not need computing $G_u$ and provides the desired accuracy without unnecessary choice of too small dead zone parameter.

\subsection{Model verification under Assumption {\mathversion {bold} $J(\theta) \le J_*$} } 
 
The Assumptions 1-4 use minimum required a priori information about the total disturbance $v$ in the framework of the robust control theory in the $\ell_1$ setting. Indeed, a priori Assumption A2 describes only the model of external disturbance and coprime factor perturbations corresponding to this theory. A priori Assumption 3 is in fact not an assumption, but formulation of control problem in the adaptive setup under no quantitative information about the external disturbance and coprime factor perturbations.  A priori Assumption 4 is a condition of robust stabilizability of the plant. It can be made as non-conservative as desired by choosing the parameter $\bar\delta$ sufficiently close to 1.   
The assumption of boundedness of the external disturbance is in the base of the theory of robust control in the $\ell_1$ setting and can not be weakened. The use of no additional information on the norm  $\delta^w$ of the external disturbance, except $\delta^w<+\infty$, shows maximum capabilities of feedback, but  has a negative consequence that the model (\ref{plant}) itself can not be falsified by data because any unacceptable dynamics of the closed loop adaptive system on any finite time interval can be explained by sufficiently large external disturbance. In any practical problem, however, too large unfalsified value of $I(\zeta_t)$, together with the inequality $I(\zeta_t)\le J(\theta)$, indicate most likely the unacceptability of the plant model under the Assumptions 1-4. That is why it is reasonable to use additionally the following a priori assumption.

\smallskip
\textbf{Assumption  A5}. \  $J(\theta) \le J_*$\,, where $J_*$ is chosen by the controller designer.
\smallskip

Assumption 5 can be used not only to falsify the model (\ref{plant}) itself under the Assumptions 1-5 after achieving the inequality $J(\zeta_t)>J_*$. Another possible application or interpretation of the Assumption 5 is to test whether the problem
\[
  \limsup_{t\to+\infty} \ |y_t| \le J_*
\]
is solvable for the model  (\ref{plant}) under the Assumptions 1-5 or, e.g., model  with larger $n$ and/or $m$ is necessary.


\section{Simulations}
\label{secsim}

In this section, we present simulations for unstable plant with the poles  0.9, 0.9, $0.8\pm 0.4i$, zeros 1.2, 1.2, and $b_1=2$, which corresponds, with $10^{-4}$ accuracy, to the coefficient vector
$
 \xi=(-4.2222, \    6.9290, \   -5.2469, \    1.5432, \    2.0000, \  -3.3333, \\    1.3889 )^T\,.
$
So the dimension of $\theta$ is 10 and the number of estimated parameters, which is the dimension of $\zeta$, equals 9. The total disturbance is modeled in the form
\begin{equation}
\label{vsim}
 v_t= w_t+0.2\delta^1_t  |y_{t-\mu}^{t-1}|+0.02\delta^2_t |u_{t-\mu}^{t-1}|\,, 
\end{equation}
where $\mu=20$, $ (\delta^w, \delta^y, \delta^u)=(1, \ 0.2,\ 0.02)$, and  $v_t$ is either \textit{random} with $w_t, \delta^1_t, \delta^2_t$ being independent and uniformly distributed on [-1,1] or $\delta^1_t=\cos(5t)$, $\delta^2_t=\sin(5t)$ in the case of \textit{deterministic perturbations}. The polytope $\Xi$ in the Assumption 1 is defined by the inequalities
\[
|a_i|\le 20,  \ i=1,2,3,4, \ |b_j|\le 10, \ j=1,2,3, \ \ b_1\ge 0.1,
\]
\[
 b_1-b_3\ge 0.01, \ b_1-b_2+b_3\ge 0.01, \ b_1+b_2+b_3\ge 0.01\,,
\]
where the inequalities in the bottom row define a compact subset of the set of stable polynomials, which is described by these inequalities with the right hand sides equal zero. The replacement of zeros by positive scalars is necessary for the compactness of a priori set $\Theta_0$ and the existence of $G_u<+\infty$ in (\ref{epsub}). The initial data $y_{-4}^{-1}$ are random, $ \xi_0=(0, \ 0, \ 0, \ 0, \ 1, \ 0, \ 0)^T$, $\varepsilon=0.001$, $\bar\mu=2\mu=40$.

Simulations for the adaptive optimal controller (\ref{ac}), (\ref{cut}) with the estimation algorithm  (\ref{zetat1})-(\ref{zetat3}) are compared with those for the adaptive controller (\ref{ac}) with the classical \textit{recursive least squares} (RLS) estimation algorithm in the form
\begin{align}
 \label{rls}
 & \xi_{t+1}=Pr_\Xi (\xi_t+K_t(y_{t+1}-\xi_t^T\varphi_t))\,, \\ \nonumber
 & K_t=\frac{P_t\varphi_t}{1+\varphi_t^T P_t\varphi}\,, \ P_{t+1}=(I-K_t\varphi_t^T)P_t\,, \ P_0=0.001 I \,,
\end{align}
where $Pr_\Xi$ denotes the projection to the nearest (under the Eucledian norm) point in $\Xi$. It is known that adaptive control based on the RLS type estimation is optimal with respect to  the mean-square type control criterion for the plant (\ref{plant}) under random external disturbance and no uncertainties \cite{guo94}. At the same time, to the best of our knowledge, there are no proven results on the stability of the RLS based adaptive control of systems under coprime factor perturbations and bounded disturbance.
In order to compare the impact of the worst-case disturbance and perturbations on the dynamics of the adaptive systems with RLS and optimal estimates, the total disturbance of the form
\begin{equation}
 \label{w-cv}
 v_{t+1}=(\delta^w+\delta^y p^y_{t+1}+\delta^u p^u_{t+1})\,\mbox{sign}(\xi_t^T\varphi_t)
\end{equation}
was modelled to \textit{maximize} next outputs $|y_{t}|$ on the time intervals [801, 810] and [1201, 1210].

\begin{figure}[h]
 {\includegraphics[height=4cm,width=9cm]{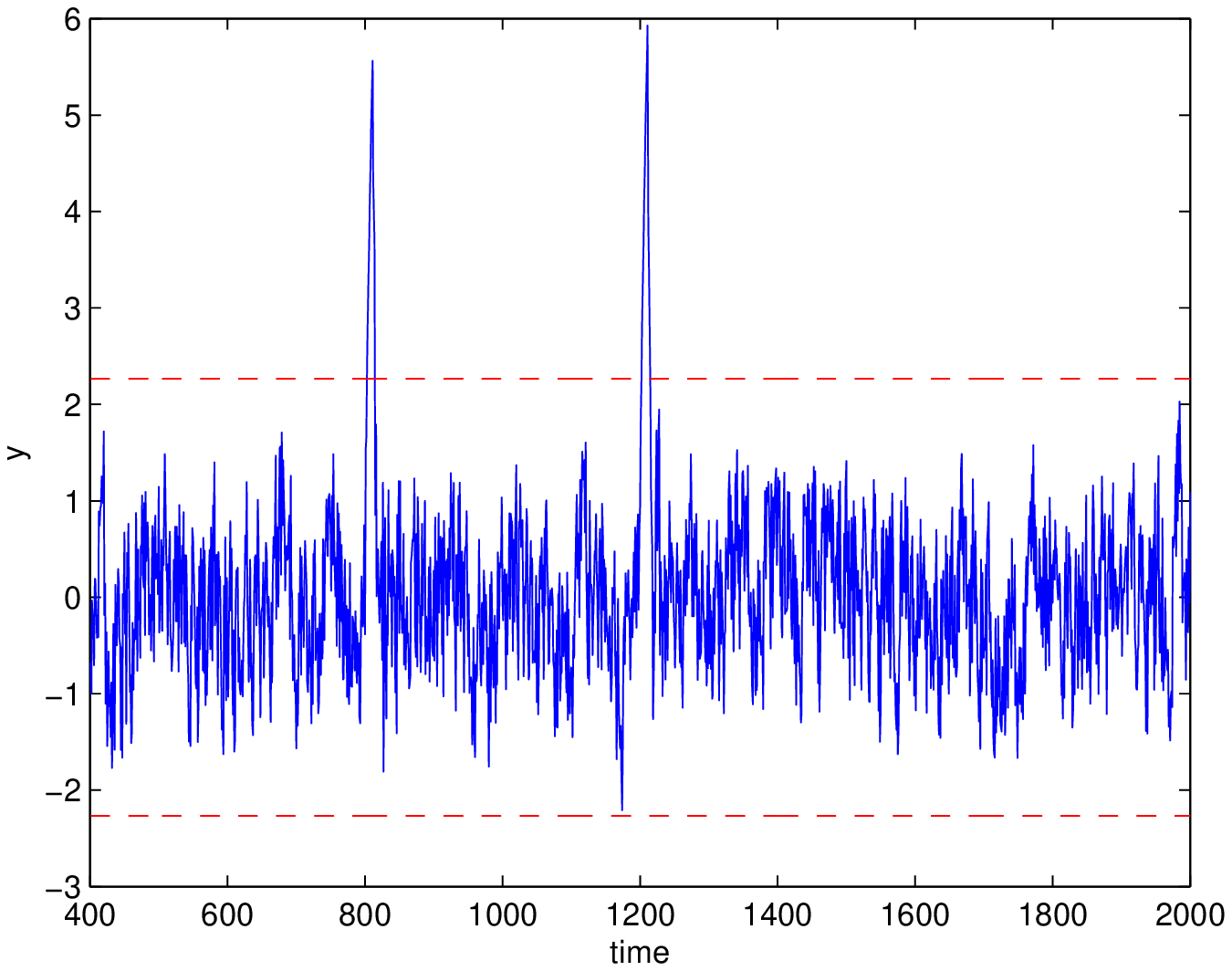}}
 {\includegraphics[height=4cm,width=9cm]{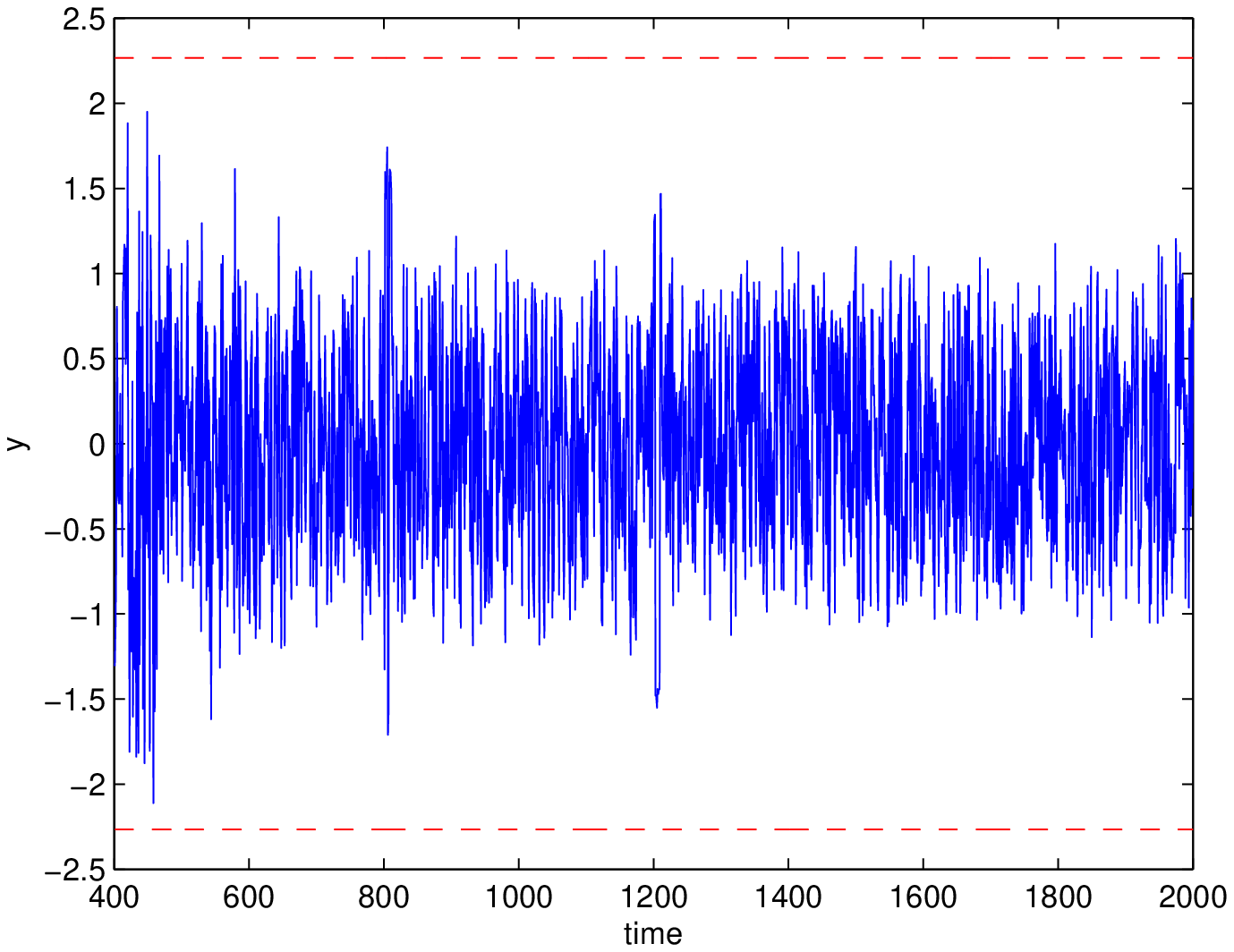}}
 \caption{Graphs of $y_t$ for the RLS algorithms (\ref{rls}) (left) and the optimal algorithm  (\ref{zetat1})--(\ref{zetat3}) (right), \ $\pm J(\theta)$ -- red dash lines.}
 \label{figrand}
\end{figure}

\begin{figure}[h]
 {\includegraphics[height=4cm,width=9cm]{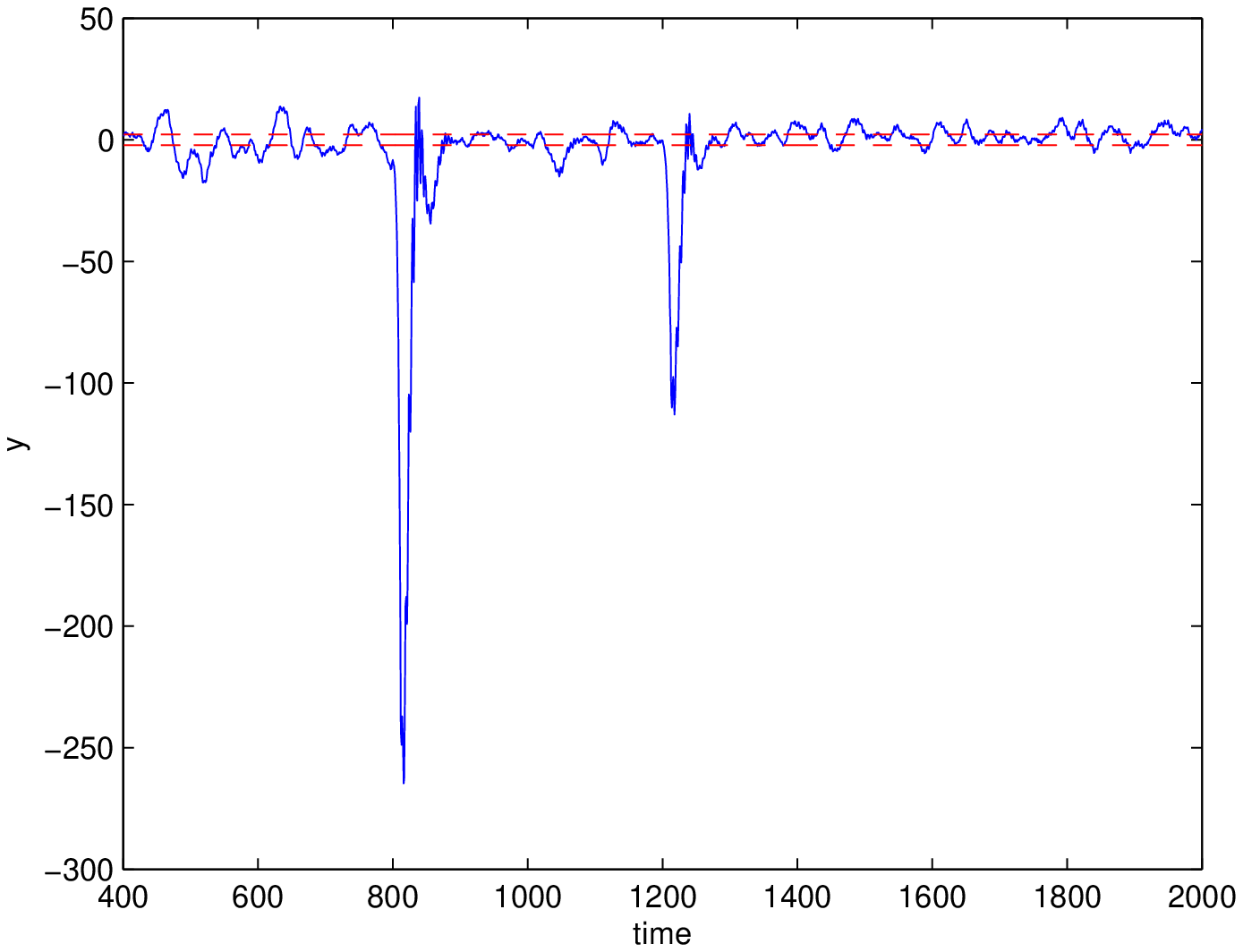}}
 {\includegraphics[height=4cm,width=9cm]{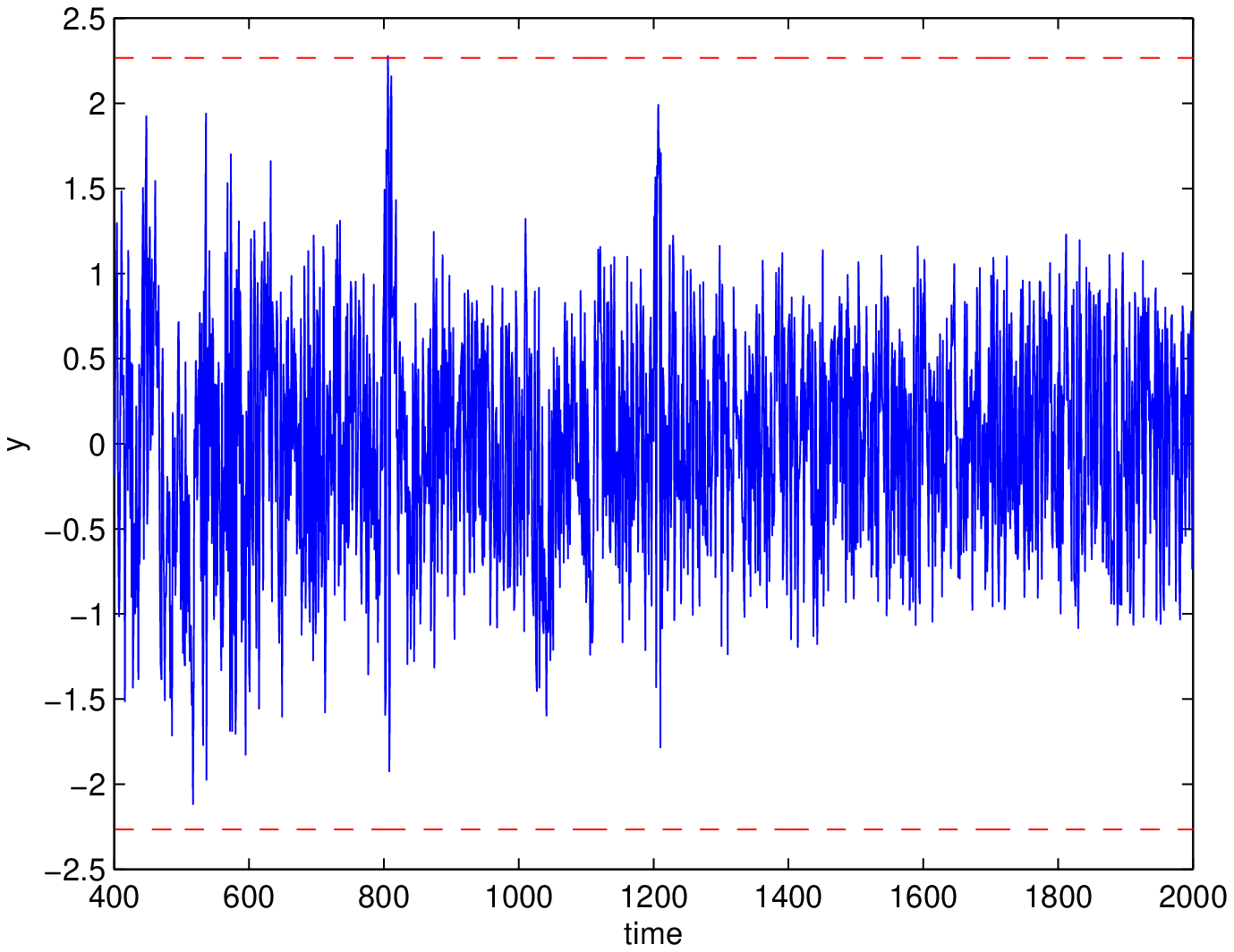}}
 \caption{Graphs of $y_t$ for the RLS algorithms (\ref{rls}) (left) and the optimal algorithm  (\ref{zetat1})--(\ref{zetat3}) (right), \ $\pm J(\theta)$ -- red dash lines.}
 \label{figdet3}
\end{figure}

Fig. \ref{figrand} presents typical graphs of the outputs $y_t$ for the adaptive systems  with the RLS estimation  (\ref{rls}) (left) and the optimal estimation (\ref{zetat1})--(\ref{zetat3}) (right) under the random total disturbance $v_t$ with the same samples $w_t, \delta^1_t, \delta^2_t$. The red dash lines on Figures \ref{figrand} and \ref{figdet3} correspond to the optimal values of the control criterion $\pm J(\theta)=\pm 2.267$. Simulations with the random total disturbance $v$ illustrate that the RLS estimation can not prevent possible bursts of the output out of the optimal interval $[-J(\theta, \ J(\theta)]$ and the second burst can be greater than the first one.

Fig. \ref{figdet3} presents graphs of the outputs $y_t$ for the adaptive systems  with the RLS estimation  (\ref{rls}) and the optimal estimation (\ref{zetat1})--(\ref{zetat3}) under the same $w_t$ and deterministic uncertainties of the form (\ref{vsim}) with $\delta^1_t=\cos(5t)$, $\delta^2_t=\sin(5t)$. In the most of simulations with the RLS estimation and random initial data $y_{t-n+1}^0$, the outputs $y_t$ in steady-state go beyond the optimal interval  $[-J(\theta),J(\theta)]$ and the bursts of $y_t$ after the worst-case total disturbances (\ref{w-cv}) remained in the interval $[-10J(\theta),10J(\theta)]$. In the specific simulation presented on Fig. \ref{figdet3}, the burst $|y_{816}|=264.695$ exceeds $J(\theta)\cdot 10^2$.

\begin{figure}[h]
 {\includegraphics[height=4cm,width=9cm]{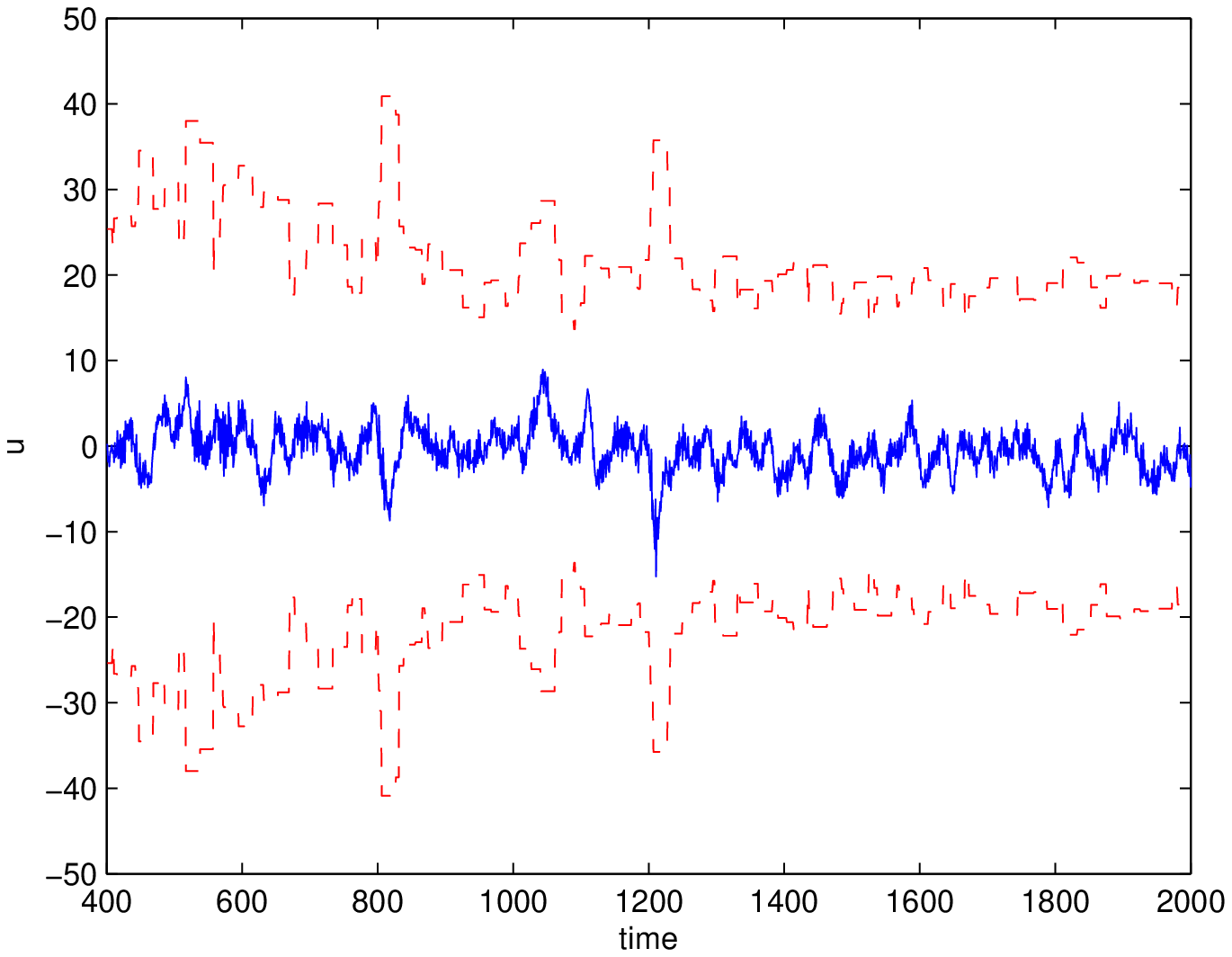}}
 {\includegraphics[height=4cm,width=9cm]{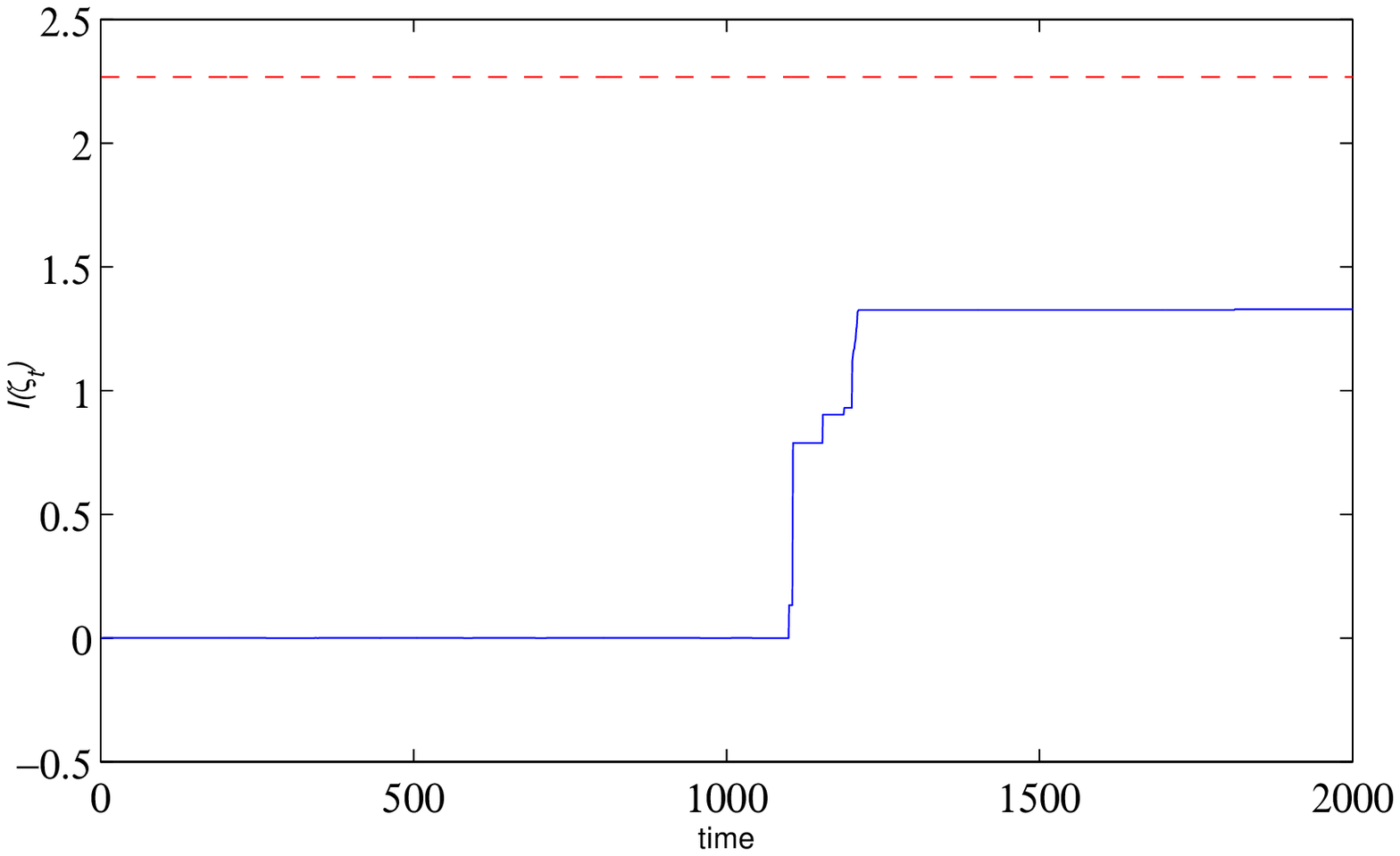}}
 \caption{Left - the graphs of $u_t$ (blue line) and $\pm\bar u_t$ (red lines); right - the graph of $I(\zeta_t)$ (blue line) and $J(\theta)$ (red line).}
 \label{figJtdet3}
\end{figure}

The left graph on Fig. \ref{figJtdet3} illustrates no violations of the inequalities (\ref{uub}). The red lines on this figure correspond to the values of $\pm \bar u_t$,  where
\[
 \bar u_t=\|G^\xi\| |y_{t+\mu-\bar\mu}^{t}|\,, \quad \mu=20\,, \quad \bar\mu=40\,.
\]
The right graph on the Fig. \ref{figJtdet3} presents the graph of the best unfalsified values of $I(\zeta_t)$ and  illustrates the model verification.  It is interesting to note, that  the final unfalsified value $I(\zeta_{2000})=1.3283$ of the control criterion is considerably less than the optimal value $J(\theta)=2.267$ despite the worst-case total disturbance of maximal magnitudes on two time intervals. The same was the case in all simulations.


Let us make some comments to simulations made on PC with the processor 8xIntel Core I7-4770 CPU @3.40 GHz. Time for simulations on the time interval [0, 2000] was typically around 0.3 sec for the RLS estimation algorithm and around 1.5 sec. for the optimal estimation. Taking time for computing the RLS estimates equal to zero, one can consider 1.2 sec as approximate time for computing the optimal estimates. The number of updates of the set estimates $Z_t$ and the vector estimates $\zeta_t$ was typically in the interval 60-70 updates and did not grow considerably with the increase of the final time of simulations up to 5000. The above figures illustrate computational tractability of the adaptive optimal control for the system with 10 unknown (and 9 estimated) parameters.

It is known that the RLS estimates for the plant (\ref{plant}) under stochastic disturbance and no perturbations converge to the line $\{ c\xi |\ \forall c\in \mathbb{R}\}$ \cite{kum85}. This is a consequence of the fact that the equation of the optimal controller (\ref{optcont}) can be rewritten in the form $\xi\phi_t=0$, which is equivalent to the equation $(c\xi)\phi_t=0\,,  \forall c\ne 0$. Proximity of the RLS and optimal estimates to this line was observed in all simulations for the plant with perturbations. For the simulations presented on Fig. \ref{figdet3}, the cosine of the angle between the vectors $\xi_{800}$ and $\xi$ equals 0.9852 for the RLS estimation and 0.9789 for the optimal estimation with
$
\xi_{800}=(-3.0147,    4.1941,   -3.0095,    1.2848,    1.8266, -1.3640,  \\  1.2172)^T
$, and
$
 \zeta_{800}=(-14.4678,   20.0,  -14.3131,    3.2475, \\   6.8781,   -9.3749,  2.9581,    0.9500,    0.0465)^T
$,
respectively. 
One can see that the $\xi_{800}$ is closer to $\xi$ and to the line $\{ c\xi | \forall c\in \mathbb{R} \}$ than  the $\xi_{800}$-component of the optimal estimate $\zeta_{800}$. However, the quality of the RLS based controller is unacceptable in practice. The reason of this difference is in the fact, that much more information in the form of polyhedral estimates is used in the computation of the optimal estimates. Note that no stochastic embedding into any estimation algorithm can  guarantee the convergence of estimates to the 'true' vector $\xi$ of the plant (\ref{plant}) under deterministic perturbations (\ref{v1}). From the deterministic robust control point of view, there are no 'true' parameters of the nominal model and the problem is to compute the best model to meet a desired control objective.


\section{Conclusion}
\label{con}

In this paper, the problem of adaptive robust optimal stabilization is considered in the optimal setting. The controlled SISO plant is described by a discrete-time linear time-invariant minimum phase nominal model under nonlinear and/or time-varying coprime factor perturbations and bounded external disturbance. The coefficients of the transfer function are assumed to be in a known polyhedron. The unknown are the coefficients of the transfer function of the nominal model, the norm of the disturbance and the gains of the coprime factor perturbations. The control criterion in the form of the worst-case steady-state upper bound on the plant output dictates the consideration of the optimal problem within the $\ell_1$-theory of robust control associated with the $\ell_\infty$ signal space and bounded disturbances. The optimal controller for the know plant depends on the coefficients of the transfer function alone, but the optimal upper bound on the output is a nonconvex function of the coefficients, the norm of additive disturbance, and the gains of perturbations. Under described a priori information, all unknown parameters are non-identifiable. Nevertheless, the proposed adaptive control guarantees, with the prescribed accuracy, the same steady-state upper bound on the plant output as the optimal controller for the known plant, that is the adaptive control realizes the maximum capability feedback with the prescribed accuracy. The solution of the optimal problem is based on the use of the control criterion as the identification criterion. Current vector estimates are computed via minimizing of the control criterion on polyhedral upper estimates  of the set of unfalsified by data parameters, the gains of perturbations and the norm of disturbance including. Computational tractability of the proposed adaptive control is illustrated by simulations.

\end{document}